\documentclass[leqno]{article}
\usepackage{amsmath,amsfonts,amsthm,indentfirst}

\newtheorem{theorem}{Theorem}
\newtheorem{lemma}[theorem]{Lemma}
\newtheorem{corollary}[theorem]{Corollary}
\newtheorem{proposition}[theorem]{Proposition}

\newcommand{\<}{\kern.0833em}
\newcommand{\adj}{\<\mathrm{adj}\kern.05em}
\newcommand{\CGr}{\mathrm{CGr}}
\newcommand{\Gr}{\mathrm{Gr}}
\newcommand{\rank}{\mathrm{rank}}
\renewcommand{\[}{\begin{equation}}
\renewcommand{\]}{\end{equation}}
\newcommand{\ph}{\varphi}
\newcommand{\bwedge}{\mbox{\large{$\wedge$}}}
\newcommand{\nxn}{$\!\<n\<{\times}\<n\!\<$}
\newcommand{\tr}{^{\mathrm T}}
\def\r#1{{\rm #1}}

\makeatletter
\let\@zzmaketitle\@maketitle
\def\@maketitle{%
  \let\@zznull\null
  \setbox\@tempboxa\hbox to\hsize{\hss\raise 20mm\hbox%
{\texttt{Comments, corrections, and related references welcomed,
as always!}}\hss}%
  \ht\@tempboxa=0pt\dp\@tempboxa=0pt%
  \def\null{\box\@tempboxa}%
  \@zzmaketitle
  \let\null\@zznull
}
\makeatother

\begin{document}

\title{Can one factor the classical adjoint of a generic
matrix?\thanks{
E-mail address of author:
gbergman@math.berkeley.edu.
This preprint is readable online at
http://math.berkeley.edu/\protect\linebreak[0]%
{$\!\sim$}gbergman\protect\linebreak[0]%
/papers/factor\_adj.\{tex,dvi\} and arXiv:math.AC/0306126\,.
(The former files will in general be updated more frequently than
the arXiv copy.)}}
\author{George M. Bergman}
\maketitle

\begin{abstract}
Let $k$ be an integral domain, $n$ a positive integer, $X$ a generic
\nxn\ matrix over $k$ (i.e., the matrix $(x_{ij})$ over a polynomial
ring $k[x_{ij}]$ in $n^2$ indeterminates $x_{ij}),$
and $\adj(X)$ its classical adjoint.
For $\mathrm{char}\;k=0$ it is shown that
if $n$ is odd, $\adj(X)$ is not the product
of two noninvertible \nxn\ matrices over $k[x_{ij}],$ while for
$n$ even, only one particular sort of factorization can occur.
Whether the corresponding result
holds in positive characteristic is open.

The operation $\adj$ on matrices arises from the $\!(n\<{-}1)\!$st
exterior power functor on modules; the analogous factorization question
is raised for matrix constructions arising from other functors.
\end{abstract}

\section{Introduction.}\label{intro}

If $A$ is an \nxn\ matrix over a commutative ring,
and $\adj(A)$ its \emph{classical adjoint},
i.e., the \nxn\ matrix of appropriately signed minors of $A,$ we
have the well-known factorization
\[\det(A)\,I_n\;=\;A\:\adj(A)\label{detAI=}\]
\mbox{\cite[p.193,~(5)]{PC.ClAl}},
\mbox{\cite[Prop.~XIII.4.16]{SL.Alg}}.
Do the factors on the right
in~(\ref{detAI=}) have any further natural factorizations?

To make this question precise, let us fix an integral domain
$k,$ and let $k[x_{ij}]$ be a polynomial ring
in $n^2$ indeterminates $x_{ij}\;(1\leq i,j\leq n).$
The matrix $X= (x_{ij})$ is
called a \emph{generic} \nxn\ matrix over $k,$
and we ask whether one can refine the factorization
\[\det(X)\,I_n\;=\;X\:\adj(X)\label{detXI=}\]
as a factorization of $\det(X)\,I_n$ into noninvertible
\nxn\ matrices over $k[x_{ij}].$

Note that the determinant of $\<\det(X)\,I_n$ is $\,\det(X)^n;$
hence, in view of~(\ref{detXI=}), $\det(\adj(X))=\det(X)^{n-1}.$
Moreover, if $k$ is an integral domain,
$\det(X)$ is irreducible over $k[x_{ij}].$
Indeed, $\det(X)$ is homogeneous of degree~$\!1\!$ in the entries of
each row of $X,$ hence any factor must be homogeneous of degree
$\!1\!$~or $\!0\!$ in those entries; hence if,
in a factorization of $\det(X),$ one factor involves
an $x$ from some row, then the other factor cannot involve any
$x$ from that row, hence the first factor must involve all the
$x$'s in that row; moreover the same applies to columns.
It follows that in any factorization, one factor
must involve all the indeterminates and the other factor none; hence
the latter belongs to $k,$ and must be a unit thereof, since
the coefficients of the monomials in $\det(X)$ are~$\pm 1.$

It follows that in~(\ref{detXI=}), $X$ cannot be factored further
into noninvertible square matrices, and that if $k$ is a
field, so that $k[x_{ij}]$ is a unique factorization domain,
any such factorization of the other term, say
\[\adj(X)\;=\;YZ, \label{ad=YZ}\]
must, up to units, satisfy
\[\det(Y)\;=\;\det(X)^d,\quad\det(Z)\;=\;\det(X)^{n-1-d},
\qquad\mbox{where $0<d<n-1.$} \label{detY,Z}\]
One can deduce that the latter statement is also true whenever
$k$ is an integral domain, by noting that it holds over
the field of fractions of $k,$ and again handling scalars by looking at
monomials over $k$ which have coefficient~$\!1\!$ in $\det(X)^n.$

Given any homomorphism $\ph$ of $k$-algebras, let us use the same
symbol $\ph$ for the induced map on \nxn\ matrices.
Note that for each \nxn\ matrix $A$ over a commutative
$\!k\!$-algebra $R,$ there is a unique $\!k\!$-algebra homomorphism
$\ph_A\!: k[x_{ij}]\rightarrow R$ carrying
the generic matrix $X$ to $A\<.$
This map $\ph_A$ will therefore carry a factorization~(\ref{ad=YZ}),
if one exists, to a factorization of $\adj(A);$ and the entries
of the factor matrices $\ph_A(Y)$ and $\ph_A(Z)$ will be given
by polynomials in the entries of $A\<.$
In particular, if $R=k=\!$ the field of real or complex numbers, the
matrices $\ph_A(Y)$ and $\ph_A(Z)$ will vary continuously with $A\<.$

By combining this observation with topological results from
\cite{CC+ZR} and \cite{G+H+S}, we shall, in Theorem~\ref{main},
exclude, for $k$ of characteristic~$0\<,$ all possible cases of
factorizations~(\ref{ad=YZ}) satisfying~(\ref{detY,Z}), except possibly
when $n$ is even and one of the
exponents $d$ or $n\<{-}1{-}\<d$ is~$1\<.$

In the first preprint version of this note, it was posed as an open
question whether the latter case actually occurred;
I and those I spoke with expected a negative answer.
However, an affirmative answer has been obtained (for arbitrary
$k)$ by R.-O.~Buchweitz and G.~Leuschke \cite{RB+GL}.
Section~\ref{even_n} below gives a quick proof of the existence of such
a factorization, inspired by working backwards from the construction
of~\cite{RB+GL}.

Let me make one caveat before beginning the development of
Theorem~\ref{main}: If we had a factorization~(\ref{ad=YZ}), the
induced factorizations $\adj(A) = \ph_A(Y)\<\ph_A(Z)$ would be
functorial, in the sense that they would respect
homomorphisms among $\!k\!$-algebras; but it cannot be assumed that they
would have other reasonable functoriality-like properties,
even when these hold for $\adj$ itself.
For instance, because the construction $\adj$ is induced by dualization
(matrix transpose) followed by the $\!(n\<{-}1)\!\<$st exterior
power functor on modules, it satisfies the multiplicative relation
\[\adj(A\<B)\ =\ \adj(B)\:\adj(A),\label{mult}\]
but it does not follow that for $Y$ as in~(\ref{ad=YZ}) we would
have $\ph_{AB}(Y) = \ph_B(Y)\,\ph_A(Y).$
For another example,~(\ref{detAI=})
applied to a matrix $U\in \mathrm{SL}_n(k),$
gives $\adj(U)=U^{-1},$ whence~(\ref{mult}) yields
$\adj(UA\,U^{-1}) = U\,\adj(A)\:U^{-1};$ but again,
no such property can be assumed for $\ph_A(Y).$

On the other hand, let us note some valid consequences of
functoriality in $k.$
If we had a factorization~(\ref{ad=YZ}) satisfying~(\ref{detY,Z}) over a
base ring $k,$ we would immediately get such a factorization over
any ring to which $k$ can be mapped homomorphically; hence in proving
{\em nonexistence} of such factorizations, results for algebraically
closed
fields $k$ will imply results for general commutative rings $k.$
Moreover, since a factorization of $\adj(X)$ over
a given $k$ involves only finitely many elements of $k,$
and any finitely generated field of characteristic~0 embeds in
$\mathbb C\<,$ restrictions on the form of
factorization with $k=\mathbb C$ will imply the
corresponding restrictions for all fields of characteristic~$\!0,\!$
and hence for all integral domains of characteristic~$0.$

(The limitation to integral domains is needed so that we can say that
any factorization~(\ref{ad=YZ}) satisfies~(\ref{detY,Z}) for some $d\<.$
Over a ring $k$ of the form $k_1\times k_2,$ in contrast, we can get a
factorization that ``looks like'' $\adj(X)=\adj(X)\cdot I_n$ over
$k_1,$ but like $\adj(X)=I_n\cdot\adj(X)$ over $k_2.$
If we consider only factorizations satisfying~(\ref{detY,Z}), on the
other hand, the restrictions on $d$ that we will obtain for
$k=\mathbb C$ will hold for any $k$ of characteristic~$0.)$

\section{Valuations and ranks.}\label{val}

We shall show below for $k$ a field of arbitrary characteristic
that if there exists a
factorization~(\ref{ad=YZ}), then, for appropriate families of
matrices $A,$ the induced matrices $\ph_A(Y)$ have constant rank.
Varying $A,$ we will get a continuous map between Grassmannian
varieties; it is
to this that we will apply topological results in the next section.
Our proof of the constant-rank result begins with

\begin{lemma}\label{DVR}
Let $R$ be a discrete valuation
ring, with valuation $v,$ maximal ideal $\mathbf m,$ and residue map
$\pi\!:\, R \rightarrow R/\mathbf m\<.$
If $M$ is an \nxn\ matrix over $R$ such that $\pi(M)$ has
nullity $r$ \r(i.e., rank $n\<{-}\<r),$ then $v(\det(M))\geq r.$
\end{lemma}\begin{proof}
Left multiplication by some invertible matrix $\pi(U)$
over $R/\mathbf m$
turns $\pi(M)$ into a matrix whose last $r$ rows are zero.
Since $\pi(U)$ is invertible, $v(\det(U))=0,$ so $v(\det(M))=
v(\det(UM)),$ which is $\geq r$ since $UM$ has $r$ rows in $\mathbf m.$
\end{proof}

\begin{corollary} \label{UFD}
Let $p$ be an irreducible element in a unique factorization domain
$R,$ $v_p$ the corresponding valuation on $R,$ and
$\pi_p: R\rightarrow R/pR$ the residue map.
Then for $M$ a square matrix over $R,$ $\pi_p(M)$
has rank at least $n- v_p(\det(M)).$
\end{corollary}\begin{proof}
Localize at $pR,$ and apply the preceding lemma in contrapositive form.
\end{proof}

Using the above results we can now prove

\begin{lemma} \label{rk}
Let $X=(x_{ij})$ be a generic \nxn\ matrix over a field $k,$ and
suppose $\adj(X)$ admits a factorization~\textnormal{(\ref{ad=YZ})}
satisfying~\textnormal{(\ref{detY,Z})} for some $d.$
Let $A$ be any \nxn\ matrix over $k$ which has the eigenvalue
$0$ with multiplicity exactly~$\!1,\!$ and let $\ph_A$ denote the
homomorphism $k[x_{ij}]\rightarrow k$ taking $X$ to $A\<.$
Then
$$\rank(\ph_A(Y))=n\<{-}\<d,\qquad\rank(\ph_A(Z))=d\<{+}1,\qquad
\rank(\ph_A(X\<Y))=n\<{-}1{-}\<d,\qquad\rank(\ph_A(ZX))=d.$$
\end{lemma}\begin{proof}
Let $k[\<t\<]$ be a polynomial ring in one indeterminate, and
$v_t$ the valuation on this ring induced by the element $t.$
From the hypothesis on $A,$ we see that $\det(tI_n + A),$ i.e., the
characteristic polynomial of $-A$ in the indeterminate $t,$ has constant
term $0$ but nonzero coefficient of $t,$ so $v_t(\det(tI_n + A))=1.$
Writing $\psi: k[x_{ij}]\rightarrow k[\<t\<]$ for
$\ph_{(tI_n+A)},$ i.e., the $\!k\!$-algebra
homomorphism taking $X$ to $tI_n+A,$ we get
\[v_t(\det(\psi(Y)))\;=\;v_t(\psi(\det(Y)))\;=\;
v_t(\psi(\det(X)^d))\;=\;v_t(\det(tI_n+A)^d)\;=\;d.\label{v_t=d}\]
Letting $\pi_t\!:\, k[\<t\<]\rightarrow k$ take $t$ to $0,$
we have $\pi_t\<\psi = \ph_A,$ hence applying
Corollary~\ref{UFD} to~(\ref{v_t=d}) we get $\rank(\ph_A(Y))\geq n-d.$
Similarly, $\rank(\ph_A(Z))\geq n-(n-1-d) = d+1.$

On the other hand, note that $A\;\ph_A(Y)\:\ph_A(Z) = \ph_A(X\,YZ)=
\ph_A(\det(X)\,I_n)=\det(A)\,I_n= 0,$
so the nullities of $A,$ of $\ph_A(Y)$ and of $\ph_A(Z)$ must add
up to at least $n,$ i.e., their ranks can sum to at most $2n.$
Since the rank of the first is $n-1$ and those of the other two are
{\em at least} $n-d$ and $d+1,$ these must be their exact values,
giving the first two equalities.
The other two are seen similarly.
(In obtaining the last one, we use~(\ref{detXI=}) in the
form $\det(X)\,I_n=\adj(X)\,X,$ easily deduced from the form given.)
\end{proof}

\textit{Remark:} The above hypothesis that the
eigenvalue~$\!0\!$ have multiplicity~$1$
is stronger than saying that $A$ has rank~$\:n-1.$
For example, the matrix consisting of a single \nxn\ Jordan block
with eigenvalue~0 has rank~$\:n-1,$ but eigenvalue~$0$ with
multiplicity~$n.$
\medskip

We will formulate our next result in algebraic-geometric terms.
For base field $\mathbb R$ or $\mathbb C\<,$ this formulation will imply
the corresponding topological statement, which is what we will actually
use in the next section; but as we will discuss in \S\ref{unproved},
the algebraic-geometric statement has the
potential of yielding results in positive characteristic as well.

For $0\leq d\leq n,$ let $\Gr_k(d,n)$ denote
the Grassmannian variety over $k$ whose $\!K\!$-valued points, for a
field $K$ over $k,$ correspond to
$\!d\!$-dimensional subspaces $V_d\subseteq K^n.$
On the other hand, let $\CGr_k(d,n)$ (for ``complemented Grassmannian'')
denote the variety whose $\!K\!$-valued points correspond to
pairs $(V_d,V'_{n-d})$ consisting of a $\!d\!$-dimensional
subspace $V_d$ and an $\!(n\<{-}\<d)\!$-dimensional subspace
$V'_{n-d}$ such that $K^n=V_d\oplus V'_{n-d}.$

The variety $\Gr_k(d,n)$ is projective; in
particular $\Gr_k(1,n)$ is $\!(n\<{-}1)\!$-dimensional projective space.
On the other hand, $\CGr_k(d,n)$ is affine, since it can be identified
with the variety of idempotent \nxn\ matrices of rank~$\:d.$

\begin{proposition}\label{CG>G}
Suppose $\adj(X)$
admits a factorization~\textnormal{(\ref{ad=YZ})}
satisfying~\textnormal{(\ref{detY,Z})} for some $d.$
Then there exists a morphism of
varieties $\CGr_k(1,n)\rightarrow\Gr_k(d,n)$ which
takes every pair $(V_1,\:V'_{n-1})$ to a subspace of its second
component~$\:V'_{n-1},$ and a morphism
$\CGr_k(1,n)\rightarrow\Gr_k(n\<{-}1{-}\<d,n)$ with the same property.
\end{proposition}\begin{proof}
Given a $\!K\!$-valued point
$a=\<(V_1,\:V'_{n-1})$ of $\CGr_k(1,n),$ let $E_a$
denote the idempotent matrix over $K$ that projects
$K^n$ onto $V'_{n-1}$ along $V_1.$
This has eigenvalue~$0$ with multiplicity~1, hence by Lemma~\ref{rk},
$E_a\;\ph_{E_a}(Y)$ has rank $n\<{-}1{-}\<d;$ so its column space,
a subspace of the column space $V'_{n-1}$ of $E_a,$ has that rank.
This construction can be seen to give a morphism of varieties,
the second of the morphisms whose existence we were to prove.

To get the first, note that taking the transpose of the equation
(\ref{ad=YZ}) and applying it to the transpose of the
matrix $X,$ we get a factorization $\adj(X) = Z'\<Y'$ with
$\det(Y')=\det(Y)$ and $\det(Z')=\det(Z).$
Applying the preceding
result to this factorization gives the desired morphism.
\end{proof}

\section{The hairy sphere raises its unkempt head.}\label{hairy}

If $n\leq 2,$ the condition $0<d<n-1$ of~(\ref{detY,Z})
cannot be satisfied, so the first case where a
factorization~(\ref{ad=YZ}) might be possible is when $n=3,\;d=1.$
Suppose we had such a factorization for $k=\mathbb R\<.$
Every point $p$ of the unit sphere $S^2$ determines a point
$(\<\mathbb R\<p,\<(\mathbb R\<p)^\perp)$ of $\CGr_{\<\mathbb R}(1,3),$
so applying to this the first morphism of Proposition~\ref{CG>G},
we would get a continuous
map $S^2\rightarrow \Gr_{\<\mathbb R}(1,3)$ that takes each $p\in S^2$
to a $\!1\!$-dimensional subspace of $(\mathbb R\<p)^\perp;$
in other words, of the tangent space to $S^2$ at $p.$
This would constitute a ``combing of a hairy sphere'', which is
known to be impossible \cite[Theorem~16.5]{MJG}, \cite[p.282]{doC},
so no such factorization exists.

(The ``hairy sphere'' result as generally formulated asserts, for even
$m,$ the nonexistence of a nowhere zero tangent vector field on $S^m.$
What the above construction would give is a map taking each
$p\in S^2$ to a point of projective $\!2\!$-space representing an
unoriented tangent direction at $p.$
But by simple connectedness of $S^2,$ we could lift this to a
map to the universal covering space of that projective plane, $S^2,$
which would determine a tangent vector field
of everywhere unit length, giving the desired contradiction.)

\section{The general result.}\label{top}

For higher $n$ and for non-real $k,$ we will use in place of the ``hairy
sphere theorem'' some results proved in \cite{CC+ZR} and \cite{G+H+S}.
As we did with $\CGr_{\<\mathbb R}(1,3)$ and $\Gr_{\<\mathbb R}(1,3)$ in
the preceding section, in the proof of the next theorem we shall regard
varieties $\Gr_{\<\mathbb C}(d,n)$ and $\CGr_{\<\mathbb C}(d,n)$ as
topological manifolds (consisting of the $\mathbb C$-valued points of
the algebraic varieties), and so be able speak of continuous
maps between them.

\begin{theorem}\label{main}
Suppose $k$ is an integral domain of characteristic~$0.$
Then if $n$ is odd, there is no
factorization~\textnormal{(\ref{ad=YZ})} of $\adj(X)$
into noninvertible matrices, while if $n$ is even, any such
factorization has one of the exponents in~\textnormal{(\ref{detY,Z})}
equal to~$1,$ i.e., has $d=1$ or $d=n-2.$
\end{theorem}\begin{proof}
As noted in the last paragraph of \S\ref{intro}, it will suffice to
prove this result for $k=\mathbb C\<.$
Let us put a Hermitian inner product on $\mathbb C\<^n;$
then $L\mapsto (L, L^\perp)$ is a continuous map
$\Gr_{\<\mathbb C}(1,n)\rightarrow \CGr_{\<\mathbb C}(1,n).$
If we have a factorization~(\ref{ad=YZ}), Proposition~\ref{CG>G}
gives a continuous map
$\CGr_{\<\mathbb C}(1,n)\rightarrow\Gr_{\<\mathbb C}(d,n)$ taking
$(V_1,\:V'_{n-1})$ to a subspace of its second component.
Composing, we get a continuous function
$\Gr_{\<\mathbb C}(1,n)\rightarrow\Gr_{\<\mathbb C}(d,n)$ taking
each $\!1\!$-dimensional subspace $L\subseteq\mathbb C\<^n$
to a $\!d\!$-dimensional subspace $L'$ of $L^\perp.$

This gives a $\!d\!$-dimensional subbundle of the tangent bundle
on $\!n\!$-dimensional complex projective space, which by
\mbox{\cite[Theorem~1.1(ii)]{G+H+S}} is possible if and only
if $n$ is even and $d=1$ or $n\,{-}\<2.$
Alternatively we may note that the map $L\mapsto L\oplus L'$ $(L'$ as
in the preceding paragraph) takes each $\!1\!$-dimensional subspace
$L$ of $\mathbb C\<^n,$ to a $\!(d\<{+}1)\!$-dimensional subspace
containing $L,$ which by \mbox{\cite[Theorem~1.5(a)]{CC+ZR}} can only
happen if $n$ is even and $d+1=2$ or $n-1,$ i.e., again $d=1$ or $n-2.$
\end{proof}

\section{The question in positive characteristic.}\label{unproved}

I do not know whether Theorem~\ref{main} remains true
if the characteristic~$0$ hypothesis is deleted.
One could hope to prove such a result by using algebraic geometry
in place of our topological arguments.

Now the analog of \cite[Theorem~1.5(a)]{CC+ZR} with morphisms
of algebraic varieties over general algebraically closed fields
in place of continuous maps of topological spaces indeed holds
[{\em ibid.,} Theorem~1.5(b)].
However, the map $L\mapsto (L, L^\perp)$ that we called on in
our proof is not a morphism
of algebraic varieties, so we cannot use it as before to connect
Proposition~\ref{CG>G} with that result.
(It is based on a Hermitian inner product, which is not bilinear but
sesquilinear; a
genuine bilinear form on $\mathbb C\<^n$ cannot be positive definite.
And if one retreats to the case $k=\mathbb R$ and tries to use a
real inner product, this will not keep its positive definiteness at
non-real points, hence will also not lead to a morphism of varieties.)
Indeed, there can be no nontrivial morphism
of algebraic varieties $\Gr_k(1,n)\rightarrow \CGr_k(1,n),$
because $\Gr_k(1,n)$ is projective while $\CGr_k(1,n)$ is affine.

What we may hope for, instead, is an analog of
\mbox{\cite[Theorem~1.5(b)]{CC+ZR}} applying directly
to morphisms $\CGr_k(1,n)\rightarrow \Gr_k(d,n).$
We remark, however, that \mbox{\cite[Theorem~1.5(b)]{CC+ZR}},
unlike \mbox{\cite[Theorem~1.5(a)]{CC+ZR}}, has only one exceptional
case for $n$ even, the case $d=n\<{-}1,$ and not the two cases
$d=2$ and $d=n\<{-}1$ as in \cite[Theorem~1.5(a)]{CC+ZR}.
Yet the example of the next section shows that both of the latter
cases occur; so the desired result would have to be weaker
than the obvious analog of \mbox{\cite[Theorem~1.5(b)]{CC+ZR}}.

\section{A factorization when $n$ is even.}\label{even_n}

We shall now see that the sorts of factorization
allowed by Theorem~\ref{main} when $n$  is even do occur.
Our argument is
inspired by the construction of Buchweitz and Leuschke \cite{RB+GL}.

\begin{lemma}\label{rkadj}
Let $R$ be a commutative integral domain, $n$ a positive integer,
and $X$ an \nxn\ matrix over $R$ having determinant~$0.$
Then
\\[2pt]
\r{(i)}~$\rank(\adj(X)) \leq 1.$
\\[2pt]
\r{(ii)}~For any alternating \nxn\ matrix $A$ over $R,$ one has
$\adj(X)\;A\;\adj(X)\tr = 0.$
\end{lemma}\begin{proof}
(i)~If $\rank(X) = n-1,$ this follows from the equation
$X\ \adj(X) = \det(X) I_n = 0.$
If $\rank(X) < n-1,$ then all minors of $X$ are zero, so $\adj(X) = 0.$

(ii)~Since $A$ is alternating, every row $r$ of $\adj(X)$ satisfies
$r A\,r\tr = 0.$
But by~(i), all rows of $\adj(X)$ are linearly dependent,
so for any two rows $r,\ r'$ of $\adj(X)$ we have $r A\,{r'}\tr = 0.$
\end{proof}

Our desired factorization is now given by part~(iii) of

\begin{theorem}\label{factor}
Let $R$ be a commutative ring, $n$ a positive integer, $X$ an
\nxn\ matrix over $R,$ and $A$ any alternating \nxn\ matrix over $R.$
Then
\\[2pt]
\r{(i)}~All entries of $\adj(X)\;A\;\adj(X)\tr$
are divisible by $\det(X).$
\\[2pt]
\r{(ii)}~$\adj(X)\;A$ is right divisible by $X\tr$ and $A\;\adj(X)$
is left divisible by $X\tr.$
\\[2pt]
\r{(iii)}~\r(Buchweitz and Leuschke {\rm\cite{RB+GL})}~If $A$ is
invertible \r(\!\<so that $n$ is necessarily even\r) then $\adj(X)$ is
right divisible by $X\tr A$ and left divisible by $A\,X\tr.$
\end{theorem}\begin{proof}
Clearly~(i) and~(ii) reduce to the case where $R$ is a polynomial ring
over the integers, $X$ a matrix of distinct indeterminates, and
$A$ an alternating matrix having distinct indeterminates for its
above-diagonal entries.
In this case, $R$ is a UFD and $\det(X)$
an irreducible element, so that $R/(\det(X))$ is an integral domain.
Applying Lemma~\ref{rkadj}(ii) to the image of the element
$\adj(X)\;A\;\adj(X)\tr$ in this domain, we get~(i).

To get~(ii), let us rewrite~(i) (still in the case where $R$ is a
polynomial ring) as
$$\adj(X)\ A\ \adj(X)\tr ~=~ Y\,(\det(X) I_n)$$
for some matrix $Y$ over $R.$
Substituting $\det(X) I_n = X\tr \adj(X)\tr$ into the right-hand side
of this equation, we can right-cancel $\adj(X)\tr$ (since it has
nonzero determinant), getting
$\adj(X)\,A = Y X\tr,$ the desired right divisibility relation.
The left divisibility statement follows by symmetry.

For $A$ invertible,~(iii) follows from~(ii)
by putting $A^{-1}$ in place of $A.$
(Note that in the resulting factorization, the factor
$X\tr A$ or $A\,X\tr$ has, up to units, determinant $\det(X)\<.$
The other factor, with determinant $\det(X)^{n-1},$ is constructed
explicitly in \cite{RB+GL}, in terms of determinantal minors.)
\end{proof}

We record the following interesting way of looking at statement~(i)
above.

\begin{corollary}\label{sup-1}
Under the general hypothesis of Theorem~\ref{factor},
if $R$ is an integral domain
and $X$ is nonsingular, then the matrix $X^{-1}A\,(X\tr)^{-1}$ over
the field of fractions of $R$ has all its entries in $\det(X)^{-1} R.$
\r(I.e., these entries, which one would {\em a priori} expect to
write using denominator $\det(X)^2,$ can in fact be
written with denominator $\det(X).)$
\end{corollary}\begin{proof}
Multiply the statement of Theorem~\ref{factor}(i) by
$(\det X)^{-2},$ recalling that by~(\ref{detXI=}),
$\det(X)^{-1}\adj(X)=X^{-1}\<,$ and hence that
$\det(X)^{-1}\adj(X\tr)=(X\tr)^{-1}\<.$
\end{proof}

Can we push the factorizations of Theorem~\ref{factor}(iii)
still further?
Suppose $A$ and $A'$ are two invertible alternating matrices over $R,$
and we write the factorizations given by that result as
\[(A X\tr)\,Y ~=~ \adj(X) ~=~ Y'\, (X\tr A').\label{AXtY=}\]
Might $Y$ and $Y'$ themselves admit nontrivial factorizations?

A look at Theorem~\ref{main} quickly eliminates all possibilities
except that $Y$ might have a right factor whose determinant (up to
units) is $\det(X)$ and/or that $Y'$ might
have a left factor with this property.
Buchweitz and Leuschke inform me, however, that they can show that
such factorizations do not occur.

Nonetheless, their result \cite[Corollary 2.4]{RB+GL} shows that in a
different sense, the two factorizations of~(\ref{AXtY=})
have a ``common refinement''; that sense being that there exist a
constant $r\in R$ and a matrix $W$ over $R$ such that
\[\adj(X) ~=~ A(rX\tr + X\tr W X\tr)A',\label{AXWXA}\]
Thus, (\ref{AXWXA}) shows both the left divisibility of $\adj(X)$
by $A X\tr,$ and its right divisibility by $X\tr A'.$

(Cf.\ \cite[Proposition~7.3(i), p.118]{PC.FRR}.
Nothing like the hypothesis of that result is satisfied here.
However that result presents a sequence of ways in which a
noncommutative ring expression can have two factorizations;
and (\ref{AXWXA}) is an instance of the $n=2$ term of that sequence.)

\section{Further questions.}\label{questions}

The factorizations of Theorem~\ref{factor}(iii) are noncanonical:  They
depend on an arbitrary invertible alternating matrix $A.$
This suggests that the context in which they would
have a natural meaning is that of a vector space given with a
nondegenerate alternating bilinear form.
It would be interesting to know whether they can in fact be given some
``functorial'' interpretation in that context.

Returning to the question with which we began this paper, but taking
a more extravagant goal than we did, we may ask whether one can describe
\emph{all} maximal factorizations of the matrix
$\det(X)\,I_n$ into noninvertible \nxn\ matrices over $k[x_{ij}].$
For any $n,$ in addition
to the factorization~(\ref{detXI=}), the same factorization with
the order of factors reversed, and the two
factorizations arising similarly from the transpose of $X,$ there is
an obvious factorization into $n$ diagonal matrices
each having determinant $\det(X){:}$
\[\det(X)\,I_n\;=\;\mathrm{diag}(\<\det(X),1,\<...\<,1)
\cdot\mathrm{diag}(1,\det(X),1,\<...\<,1)\<\cdot\;\ldots\;
\cdot\<\mathrm{diag}(1,\<...\<,1,\det(X)\<).\label{0+0+}\]

When $n$ is odd, are~(\ref{detXI=}), its three variants noted above,
and~(\ref{0+0+}), ``essentially'' all there are?

A factorization can be trivially perturbed by multiplying
any two successive factors on the right and the left respectively
by an invertible matrix $U$ over $k[x_{ij}]$ and its inverse.
Also, because $\det(X)\,I_n$ is central, we can left-multiply the first
factor in a factorization by such a matrix $U,$ and right-multiply
the last factor by $U^{-1}.$
So we may ask whether the five factorizations we have described
form a set of representatives of the orbits of all maximal
factorizations of $\det(X)I_n$ under these sorts of perturbations.
For $n$ even, we get additional factorizations from
Theorem~\ref{factor}(iii), this time parametrized by
an alternating matrix $A.$
(It is not clear how many degrees of freedom these additional
families have, modulo the equivalence relation introduced above.)

In each of the explicit factorizations noted above, one can see or show
by homogeneity arguments that the degree of
the product matrix $\det(X)\,I_n$ in the $n^2$
indeterminates is precisely the maximum of
the sums of the degrees of the matrix entries that get multiplied
together; i.e., that there is not too much ``cancellation'' in the
calculation of $\det(X)\,I_n$ as a product.
We can, however, easily destroy this property by interpolating
invertible matrices $U$
over $k[x_{ij}]$ with entries of high degree, and their inverses.
Might there, nonetheless, be some principle saying that any
factorization of a ``good'' matrix over a polynomial ring is a
perturbation, via interpolated matrices and their
inverses, of a factorization in which the degree is well-behaved?

Turning in a different direction, let us observe that for an
\nxn\ matrix $A$ over a commutative ring $k,$ say representing a
linear map $a\!:\< k^n\rightarrow k^n,$
the classical adjoint $\adj(A)$ can be characterized as the
transpose of the matrix representing the linear map
$\bwedge\!^{n-1}\,a\!:\,
\bwedge\!^{n-1}\,k^n \rightarrow \bwedge\!^{n-1}\,k^n,$ where
$\bwedge\!^{n-1}$ denotes the $\!(n\<{-}1)\!$st exterior power functor.
If instead we apply to $a$ a lower exterior power functor
$\bwedge\!^m,$ we get an endomorphism of the
module $\bwedge\!^m\,k^n,$ which is free of rank $\binom n{\<m\<}.$
Again taking for $A$ a generic matrix $X,$ we may ask whether the
resulting $\binom n{\<m\<}{\times}\binom n{\<m\<}$ matrix over
$k[x_{ij}]$ can be factored into noninvertible square matrices.
(This matrix, incidentally, has determinant
$\det(X)\!^{\textstyle\binom {n-1}{\<m-1\<}}\!,$
and its product with the transpose of the
matrix representing $\bwedge\!^{n-m}\,a,$ with rows and columns
appropriately indexed, is $\det(X)$ times the
$\binom n{\<m\<}{\times}\binom n{\<m\<}$ identity matrix.)

Each of the above functors $\bwedge\!^m$ is a subfunctor
of the $\!m\!$-fold tensor product functor
$\otimes^m.$
Indeed, when $\mathrm{char}\;k=0,$ $\otimes^m$ decomposes
into a direct sum of subfunctors indexed by Young diagrams; the functor
$\bwedge\!^m$ corresponds to the height-$\!m\!$ column of boxes.
(The length-$\!m\!$ row of boxes likewise corresponds to the
$\!m\!$th symmetric power functor.)
We may thus pose for any such subfunctor of $\otimes^m$
the same question we have studied here for~$\bwedge\!^{n-1}\,!$

\section{Acknowledgements}

I am indebted to Zinovy Reichstein
for an extensive and valuable correspondence on this subject,
and to Ragnar Buchweitz and Graham Leuschke for their surprising
answer to the main question raised in the first version of this note.


\begin{thebibliography}{1}

\bibitem{RB+GL} Ragnar-Olaf Buchweitz and Graham Leuschke,
\textit{The adjoint of an even size matrix factors}, 7pp.,
March~23, 2004,
//www.math.toronto.edu/gleuschk/AdjointFactor.dvi\,.

\bibitem{PC.ClAl} P. M. Cohn,
\textit{Classic Algebra}, J. Wiley \& Sons, 2000.
(Revision of \textit{Algebra}, v.1, MR~\textbf{83e}:00002.)

\bibitem{PC.FRR} P. M. Cohn,
\textit{Free Rings and Their Relations},
2nd ed., London Math. Soc. Monographs, vol.\,19, Academic Press, 1985.
~MR~\textbf{87e}:16006.
(To locate in the forthcoming 3rd edition the result cited above,
look for ``leapfrog construction'' in the index, and go to the first
display in the first proposition after that construction is introduced.)

\bibitem{CC+ZR} C.~de~Concini and Z.~Reichstein,
\textit{Nesting maps of Grassmannians},
to appear, Rendiconti di Matematica, Accademia dei Lincei.
Preprint, 10\,pp., revised Nov.~2003, at
http://\linebreak[0]www.math.ubc.ca/\linebreak[0]%
$\!\sim\!$reichst/\linebreak[0]pub.html\,.

\bibitem{doC} Manfredo P. do Carmo,
\textit{Differential Geometry of Curves and Surfaces},
Prentice-Hall, 1976.
~MR~\textbf{52}\#\linebreak[0]15253.

\bibitem{G+H+S} H.~Glover, W.~Homer and R.~Stong,
\textit{Splitting the tangent bundle of projective space},
Indiana University Math Journal,
\textbf{31} (1982), 161--166.
~MR~\textbf{83f}:57016.


\bibitem{MJG} Marvin~J.~Greenberg and John~R.~Harper, 
\textit{Algebraic topology, A first course,}
Mathematics Lecture Note Series, v.58,
Benjamin/Cummings, 1981.
~MR~\textbf{83f}:55001

\bibitem{SL.Alg} Serge Lang,
\textit{Algebra},
Addison-Wesley, third edition, 1993,
reprinted as Springer Graduate Texts in Mathematics, v.211, 2002.
~MR~\textbf{2003e}:00003.

\end{thebibliography}
\end{document}